\documentclass[journal,10pt]{IEEEtran}
\usepackage[top=1in,bottom=1.05in,right=0.75in,left=0.75in]{geometry}
\makeatletter
\def\ps@headings{%
\def\@oddhead{\mbox{}\scriptsize\rightmark \hfil \thepage}%
\def\@ adversarynhead{\scriptsize\thepage \hfil \leftmark\mbox{}}%
\def\@oddfoot{}%
\def\@ adversarynfoot{}}
\makeatother
\ifCLASSINFOpdf
\else
 
\fi
 
\usepackage{epsfig}
\usepackage{array}
\newcolumntype{L}[1]{>{\raggedright\let\newline\\\arraybackslash\hspace{0pt}}m{#1}}
\newcolumntype{C}[1]{>{\centering\let\newline\\\arraybackslash\hspace{0pt}}m{#1}}
\newcolumntype{R}[1]{>{\raggedleft\let\newline\\\arraybackslash\hspace{0pt}}m{#1}}
\usepackage{amsmath}
\usepackage[short]{optidef}
\usepackage{amsthm} 
\usepackage{mathrsfs}
\newcommand{\bc}{\begin{center}}
\newcommand{\ec}{\end{center}}
\newcommand{\be}{\begin{equation}}
\newcommand{\ee}{\end{equation}}
\usepackage{flexisym}
\usepackage{algorithmicx}
\usepackage{algorithm}
\usepackage{booktabs}
\usepackage{comment}
\usepackage{soul}
\usepackage{algpseudocode}

\usepackage{amsmath}

\newcommand{\bnu}{\begin{enumerate}}
\newcommand{\enu}{\end{enumerate}}
\usepackage{cite}
\usepackage{url}
\usepackage{breqn}
\usepackage{lipsum}
\usepackage{mathtools}
\usepackage{amsmath}
\usepackage{caption}
\usepackage{subcaption}
\usepackage{algpseudocode}
\usepackage{soul}
\usepackage{blindtext, graphicx}
\usepackage{cuted}
\usepackage{color}
\usepackage{algpseudocode}

\usepackage{amsfonts}
\usepackage{lipsum}
\usepackage{cuted}
\newtheoremstyle{case}{}{}{}{}{}{:}{ }{}
\begin{document}
\pagestyle{empty}
\newgeometry{top=1in,bottom=1.05in,right=0.75in,left=0.75in}
\title{Quantum Noise-Aware RIS-Aided Wireless Networks Using Variational Encoding and Signal Stabilization}
\author{Shakil Ahmed,~\IEEEmembership{Member,~IEEE}
\vspace*{-1.05 cm}
\thanks{ Shakil Ahmed is with the Department of Electrical and Computer Engineering, Iowa State University, Ames, Iowa, USA. (email: \{shakil\}@iastate.edu).
}}

 \markboth{23rd IEEE Consumer Communications and Networking Conference (CCNC) (accepted)}
{}
\maketitle

\IEEEpeerreviewmaketitle

\begin{abstract}
This paper presents a noise-aware quantum-assisted framework for blockage prediction in reconfigurable intelligent surface (RIS)-enabled wireless networks. The proposed architecture integrates a Quantum Base Station (QBS), a Quantum RIS (QRIS), and a mobile Quantum User Node (QUN). Visual information captured by an onboard RGB camera is amplitude-encoded into quantum states, while channel state observations are mapped into quantum rotation-encoded features. These hybrid inputs are processed through variational quantum circuits, enabling ternary classification of the link status. To address the inherent imperfections of noisy intermediate-scale quantum (NISQ) hardware, the system explicitly models depolarizing and dephasing channels along direct and QRIS-assisted paths. A fidelity-aware training objective is employed to jointly minimize classification loss and quantum state degradation, with amplitude damping and synthetic noise injection enhancing robustness. Simulation results on a quantum-adapted version of the ViWi dataset demonstrate that the proposed hybrid quantum model achieves superior accuracy and stability under realistic noise conditions, outperforming baseline and single-modality approaches.
\end{abstract}

\section{Introduction}

Reconfigurable Intelligent Surfaces (RIS) have emerged as a promising paradigm for enhancing wireless connectivity by enabling intelligent control over electromagnetic propagation \cite{di2020smart,huang2019reconfigurable}. Classical RIS frameworks have demonstrated the ability to redirect signals and mitigate blockage through passive beam steering and environment-aware optimization \cite{wu2020intelligent}. However, such approaches remain fundamentally constrained by classical signal processing limits, especially in highly dynamic or blocked environments. Recent advances in quantum machine learning and quantum-enhanced sensing provide new opportunities to overcome these limitations by leveraging quantum systems' representational power and potential robustness \cite{schuld2019quantum,biamonte2017quantum}.

This paper proposes a novel quantum-assisted blockage prediction framework that integrates RIS control with quantum encoding and inference. The system comprises a Quantum Base Station (QBS), a Quantum Reconfigurable Intelligent Surface (QRIS), and a mobile Quantum User Node (QUN). The QBS collects visual data using an RGB camera and measures RIS-assisted signal performance. These observations are encoded into quantum states using a hybrid format: image features are amplitude-encoded, and channel measurements are mapped via rotation gates. The resulting hybrid quantum states are processed by a variational quantum circuit (VQC), which performs ternary classification of link status—blocked, unblocked, or absent \cite{schuld2020circuit,mitarai2018quantum}.
Unlike prior models, we explicitly model quantum noise using depolarizing and dephasing channels to reflect realistic NISQ conditions \cite{preskill2018quantum}. To mitigate the effects of noise, we introduce amplitude damping during input encoding and a fidelity-aware loss function that penalizes divergence from the ideal state trajectory. This noise-regularized loss improves stability and prevents overfitting to noisy training distributions. Training is performed via the parameter-shift rule with projected updates to enforce fidelity and damping constraints \cite{crooks2019gradients}.

The authors in \cite{di2020smart} and \cite{huang2019reconfigurable} have shown that RIS can significantly enhance wireless communication by enabling controllable signal reflection in blocked or dynamic environments. Machine learning methods have been integrated into RIS systems for beam selection and blockage prediction, including the use of visual input datasets such as ViWi \cite{alrabeiah2020viwi} and learning-based optimization frameworks \cite{elbir2022deep, 10017894}. Recently, researchers have explored quantum machine learning (QML) as a new paradigm for representing high-dimensional data and solving signal inference problems more efficiently \cite{schuld2019quantum, biamonte2017quantum}. Quantum-enhanced optimization has also been applied to RIS configuration in complex propagation scenarios \cite{colella2024quantum}, including in modular settings for STAR-RIS \cite{narottama2023modular}. However, most of these approaches either assume ideal quantum channels or do not explicitly model the effect of quantum noise. Studies such as \cite{preskill2018quantum} and \cite{crooks2019gradients} highlight the critical importance of accounting for noise in near-term quantum devices. Hybrid encoding strategies for combining classical visual and channel information into quantum circuits have recently been proposed in \cite{alhassoun2024beyond}. However, noise-resilient RIS-based inference under NISQ constraints remains largely unaddressed. 

Experimental results on a quantum-adapted version of the ViWi dataset \cite{alrabeiah2020viwi} show that the proposed hybrid quantum model achieves superior accuracy and fidelity under noise, outperforming classical and single-modality baselines. These results demonstrate the feasibility and robustness of combining quantum encoding with RIS routing for predictive wireless intelligence in realistic quantum environments.
The main contributions of this work are summarized as follows:
\begin{itemize}
    \item We propose a novel quantum-assisted blockage prediction framework that integrates a QBS, a QRIS, and a mobile QUN, enabling intelligent wireless link inference under dynamic propagation conditions.

    \item A hybrid quantum encoding strategy is developed, where environmental images are amplitude-encoded and RIS-assisted signal rates are mapped to rotation-based quantum states. This results in a composite input that captures both spatial and channel-domain information.

    \item To address quantum hardware limitations, we explicitly model quantum noise through depolarizing and dephasing channels in direct and QRIS-assisted paths. We incorporate amplitude damping in the encoding and introduce a fidelity-aware loss function to enhance training robustness.

    \item Extensive simulations on a quantum-adapted version of the ViWi dataset demonstrate the effectiveness of the proposed approach. The hybrid noise-aware model outperforms classical baselines in prediction accuracy and quantum fidelity under realistic noise conditions.
\end{itemize}

\section{System Model}
\begin{figure}[!h]
\centering
\includegraphics[width=2.70in]{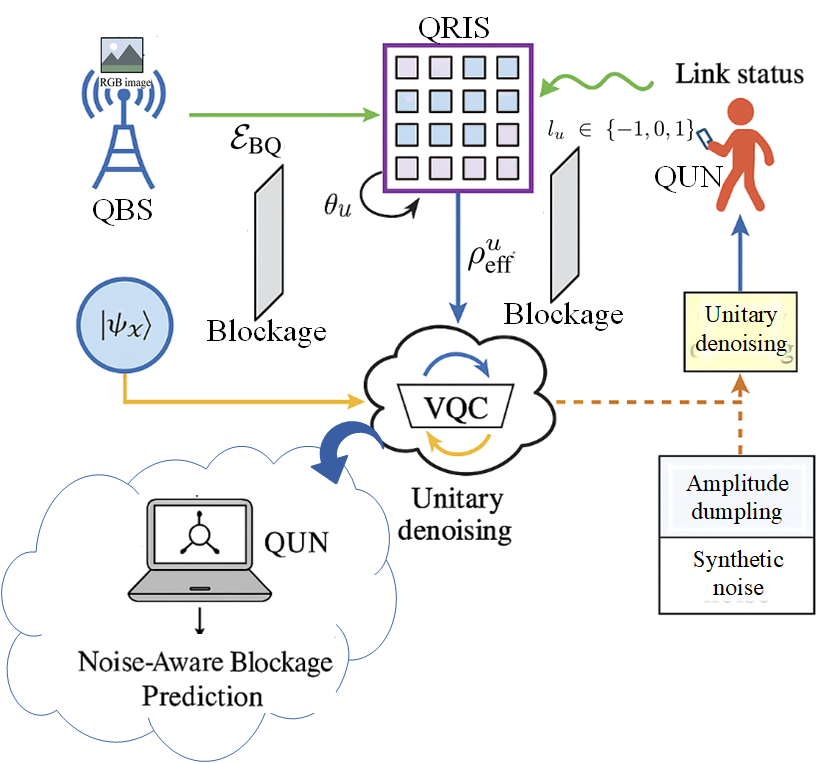}
\caption{System model of the proposed quantum RIS-assisted blockage prediction framework. The QBS encodes visual input $|\psi_X\rangle$ and channel parameter $\theta_u$ into a hybrid quantum state. The signal propagates via a direct channel $\mathcal{E}_{\text{BQ}}$ or is reflected through the QRIS, resulting in an effective received state $\rho_{\text{eff}}^u$. The VQC receives both quantum input and noisy channel state, and performs ternary blockage prediction $l_u \in \{-1, 0, 1\}$. }
\label{System_Model2}
\end{figure}
We consider a quantum-assisted wireless network composed of a Quantum Base Station (QBS), a Quantum RIS (QRIS), and a mobile Quantum User Node (QUN), as shown in Fig.~1. The QBS is responsible for both quantum signal transmission and environmental visual sensing. The QRIS functions as a programmable optical interface that applies unitary transformations to quantum states in transit, enabling intelligent signal redirection in the presence of physical obstacles \cite{di2020smart,huang2019reconfigurable}. A mobile QUN receives these transformed quantum states and executes measurement-based classification to infer the status of the wireless link. This model extends the quantum RIS framework by incorporating noise at the channel and gate levels. We propose a hybrid signal-visual encoding mechanism robust to quantum noise \cite{preskill2018quantum,schuld2019quantum}. The system architecture is illustrated in Fig.~1. Table~\ref{tab:symbols} summarizes the key quantum states and system parameters used throughout the paper.

\textit{Quantum RIS-Aided Communication Scenario:}
The QBS emits a photonic quantum signal $\rho_0$, which may propagate directly to the QUN or be reflected via the QRIS. In the direct path, the quantum signal travels through a line-of-sight (LoS) quantum channel $\mathcal{E}_{\text{BQ}}$, subject to environmental noise including scattering, decoherence, and phase errors. In the indirect path, the quantum signal first reaches the QRIS, which applies a deterministic unitary operation $U_{\text{QRIS}}$. It is then transmitted through a second noisy channel $\mathcal{E}_{\text{RQ}}$ toward the QUN. The received quantum state $\rho_{\text{eff}}^u$ at the QUN is modeled as a convex mixture of these two channels, governed by a location-dependent weighting parameter $\alpha \in [0,1]$:
\begin{equation}
\label{Eq_QS}
\rho_{\text{eff}}^u = \alpha  \mathcal{E}_{\text{BQ}}(\rho_0) + (1 - \alpha)  \mathcal{E}_{\text{RQ}}(U_{\text{QRIS}} \rho_0 U_{\text{QRIS}}^\dagger).
\end{equation}

Eq.~\ref{Eq_QS} captures the joint contribution of the direct and QRIS-assisted paths while accounting for the transformation imposed by the intelligent surface.
\begin{table}[t]
\centering
\caption{List of Symbols Used in the Paper}
\label{tab:symbols}
\begin{tabular}{|l|l|}
\hline
\multicolumn{2}{|c|}{\textbf{Quantum States}} \\
\hline
$|\psi_X\rangle$ & Amplitude-encoded quantum image state \\
$|\tilde{\psi}_X\rangle$ & Denoised quantum image input (scaled by $\gamma$) \\
$|\phi_u\rangle$ & Channel-encoded quantum state (via $R_y(\theta_u)$) \\
$|\Psi_u\rangle$ & Full hybrid input state: $|\tilde{\psi}_X\rangle \otimes |\phi_u\rangle$ \\
$\rho_0$ & Initial state emitted by the QBS \\
$\rho_{\text{eff}}^u$ & Effective received quantum state  \\
$\rho_{\text{ideal}}$ & Ideal (noise-free) output quantum state \\
$\rho_{\text{noisy}}$ & Noisy output state after propagation \\
\hline
\multicolumn{2}{|c|}{\textbf{System Parameters and Operators}} \\
\hline
$\Gamma$ & Parameters of VQC \\
$\gamma$ & Amplitude damping coefficient for input state \\
$\theta_u$ & Rotation angle for RIS-assisted rate encoding \\
$R_u[k]$ & Observed data rate for user $u$ at index $k$ \\
$Q_\Gamma(\cdot)$ & Output distribution of the VQC \\
$\mathcal{E}_{\text{BQ}}$, $\mathcal{E}_{\text{RQ}}$ & Quantum channels: QBS→QUN and QRIS→QUN \\
$\mathcal{E}_{\text{dep}}$, $\mathcal{E}_{\text{phase}}$ & Depolarizing and dephasing noise maps \\
$p$, $q$ & Depolarization and dephasing noise rates \\
$F(\cdot,\cdot)$ & Quantum fidelity between two density matrices \\
$\lambda$ & Weight of the fidelity penalty in the loss \\
$\alpha$ & Weighting factor for direct and QRIS signals \\
$l_u$ & Ternary link label: $\{-1, 0, 1\}$ \\
$H$, $W$, $C$ & Height, width, and color channels of the image \\
$N$ & Total image vector components ($N = HWC$) \\
\hline
\end{tabular}
\end{table}

\textit{Quantum Visual Encoding:}
The QBS has an RGB camera that captures a visual snapshot $X \in \mathbb{R}^{H \times W \times C}$ of the propagation environment. The image is first flattened and normalized to a unit vector $\hat{\mathbf{x}} \in \mathbb{R}^N$, with $N = H  W  C$. This vector is then amplitude-encoded into a quantum state using:
\begin{equation}
|\psi_X\rangle = \sum_{i=1}^{N} \hat{x}_i |i\rangle.
\end{equation}
Due to inherent imperfections in quantum gate operations and the increased sensitivity of amplitude encoding to over-rotation errors, the normalized image vector is scaled by a denoising factor $\gamma \in (0,1]$ to reduce circuit-level noise. The resulting stabilized quantum visual input is expressed as
\begin{equation}
|\tilde{\psi}_X\rangle = \sum_{i=1}^{N} \gamma \hat{x}_i |i\rangle,
\end{equation}
where $\hat{x}_i$ denotes the $i$-th normalized amplitude component. 

\textit{QRIS-Assisted Channel State Encoding:}
The QRIS-assisted data rate $R_u[k]$ observed at the QBS for user $u$ at time $k$ is normalized into a phase rotation angle $\theta_u \in [0, \pi]$. This angle is used to construct a single-qubit quantum state using a $R_y(\theta_u)$ rotation gate:
\begin{equation}
|\phi_u\rangle = R_y(\theta_u)|0\rangle.
\end{equation}
The resulting state captures link strength information in a format suitable for quantum classification. The complete hybrid quantum input is formed by taking the tensor product, which jointly encodes environmental and channel-level context as follows:
$|\Psi_u\rangle = |\tilde{\psi}_X\rangle \otimes |\phi_u\rangle.$

\textit{Quantum User Node and Blockage Classification:}
The QUN receives both the hybrid visual-channel quantum input $|\Psi_u\rangle$ and the effective quantum state $\rho_{\text{eff}}^u$ through the noisy propagation medium. A VQC, parameterized by trainable unitaries $\Gamma$, is applied to the quantum input state. The VQC performs layered transformations and partial measurements to predict the link status $l_u \in \{-1, 0, 1\}$, corresponding to absent, blocked, or unblocked conditions. The classifier is optimized using a hybrid quantum-classical training loop with noise-injected inputs during training, allowing it to generalize under realistic quantum hardware constraints.

\textit{Quantum Noise Modeling and Minimization:}
The direct and reflected transmission paths are modeled using noisy quantum channels to characterize the non-ideal behavior of quantum communication and processing. Specifically, the quantum links $\mathcal{E}_{\text{BQ}}$ (from the QBS to the QUN) and $\mathcal{E}_{\text{RQ}}$ (from the QRIS to the QUN) are each constructed as a serial composition of depolarizing and dephasing channels. The depolarizing map is defined as
\begin{equation}
\mathcal{E}_{\text{dep}}(\rho) = (1 - p)\rho + \frac{p}{3}(X\rho X + Y\rho Y + Z\rho Z),
\end{equation}
where $p \in [0,1]$ denotes the depolarization probability. The dephasing channel is given by
\begin{equation}
\mathcal{E}_{\text{phase}}(\rho) = (1 - q)\rho + q Z\rho Z,
\end{equation}
where $q$ controls the magnitude of phase noise. The overall channel model is expressed as $\mathcal{E} = \mathcal{E}_{\text{phase}} \circ \mathcal{E}_{\text{dep}}$ and is applied independently to both propagation links based on the user location and signal routing path.
The proposed architecture incorporates three complementary strategies to mitigate the effects of such quantum noise. First, a denoising coefficient is introduced during quantum visual state preparation, scaling the input amplitude vector to suppress gate over-rotation and circuit-induced noise. Second, synthetic quantum noise is injected into the training pipeline via simulated channel perturbations, enabling the classifier to generalize under stochastic quantum conditions. Third, the variational quantum classifier is optimized using a hybrid loss that includes both a cross-entropy term for classification and a fidelity penalty

\section{Blockage Prediction}
This work aims to predict the link status between the QBS and the QUN under dynamic environmental and channel conditions while accounting for quantum noise that arises during signal encoding, propagation, and classification. The prediction task is cast as a ternary classification problem, where the output label $l_u \in \{-1, 0, 1\}$ denotes an absent, blocked, or unblocked link state for user $u$, respectively. Given a hybrid quantum input $|\Psi_u\rangle$ composed of the denoised visual state $|\tilde{\psi}_X\rangle$ and the channel-encoded state $|\phi_u\rangle$, and an effective quantum signal state $\rho_{\text{eff}}^u$, the objective is to optimize the VQC such that the predicted class probabilities closely match the true link status.
To formalize this, let $Q_\Gamma(|\Psi_u\rangle)$ represent the output probability distribution of the VQC parameterized by $\Gamma$, and let $b_i$ denote the predicted probability for class $i$. The classification loss is defined using the categorical cross-entropy:
\begin{equation}
\mathcal{L}_{\text{CE}} = - \sum_{i=1}^{3} y_i \log b_i,
\end{equation}
where $y_i$ is a one-hot encoding of the true label $l_u$. However, minimizing only this loss in the presence of quantum noise is insufficient to guarantee stable inference. Therefore, we incorporate a fidelity-based penalty term into the objective to account for the divergence between the noisy and ideal quantum state evolutions:
\begin{equation}
\mathcal{L}_{\text{F}} = 1 - F(\rho_{\text{ideal}}, \rho_{\text{noisy}}),
\end{equation}
where $F(\cdot, \cdot)$ denotes the quantum state fidelity. The combined loss becomes:
\begin{equation}
\mathcal{L}_{\text{total}} = \mathcal{L}_{\text{CE}} + \lambda \mathcal{L}_{\text{F}},
\end{equation}
with $\lambda > 0$ controlling the influence of the fidelity regularization.
To reflect the physical and operational constraints of a noise-aware quantum system, we impose the following conditions:

 \textit{Input Damping Constraint:} The denoising factor $\gamma$ applied to the visual input must remain within a stability-preserving range to ensure amplitude normalization and avoid quantum gate saturation:
    \begin{equation}
    0 < \gamma \leq \gamma_{\max} < 1.
    \end{equation}
    
\textit{Noise Budget Constraint:} The overall effective noise level $\bar{p}$ encountered across the propagation and circuit execution paths must not exceed a threshold that would degrade the classifier's fidelity below an acceptable limit $F_{\min}$:
    \begin{equation}
    F(\rho_{\text{ideal}}, \rho_{\text{noisy}}) \geq F_{\min}.
    \end{equation}
\begin{algorithm}[h]
\caption{Noise-Aware Quantum Training for Blockage Prediction}
\begin{algorithmic}[1]
\State \textbf{Input:} Training data $\mathcal{D} = \{(X_u, R_u[k], l_u)\}$, noise rates $p$, $q$, fidelity threshold $F_{\min}$, learning rate $\eta$, penalty weight $\lambda$, damping max $\gamma_{\max}$
\State \textbf{Output:} Optimized VQC parameters $\Gamma$ and damping coefficient $\gamma$
\State Initialize VQC parameters $\Gamma$ and $\gamma \in (0, \gamma_{\max}]$
\For{each epoch}
    \For{each sample $(X_u, R_u[k], l_u)$ in $\mathcal{D}$}
        \State Flatten and normalize image $X_u$ to $\hat{\mathbf{x}}$
        \State Apply damping: $\tilde{\mathbf{x}} \gets \gamma \cdot \hat{\mathbf{x}}$
        \State $|\tilde{\psi}_X\rangle \gets$ AmplitudeEncode($\tilde{\mathbf{x}}$)
        \State Encode channel state $|\phi_u\rangle \gets R_y(\theta_u)|0\rangle$
        \State Form hybrid input: $|\Psi_u\rangle = |\tilde{\psi}_X\rangle \otimes |\phi_u\rangle$
        \State Simulate quantum noise: $\rho_{\text{noisy}} \gets \mathcal{E}(|\Psi_u\rangle)$
        \State Evaluate circuit: $b_i \gets Q_\Gamma(|\Psi_u\rangle)$
        \State Compute total loss:
        \[
        \mathcal{L}_{\text{total}} = \mathcal{L}_{\text{CE}} + \lambda(1 - F(\rho_{\text{ideal}}, \rho_{\text{noisy}}))
        \]
        \If{$F < F_{\min}$}
            \State Increase penalty: $\lambda \gets \lambda \cdot 1.1$
        \EndIf
        \State Compute gradients via parameter-shift and finite-difference
        \State Update parameters: $\Gamma \gets \Gamma - \eta \nabla_\Gamma$, \quad $\gamma \gets \text{clip}(\gamma)$
    \EndFor
\EndFor
\State \Return $\Gamma$, $\gamma$
\end{algorithmic}
\end{algorithm}

The joint training objective of the system is defined as a constrained optimization problem that minimizes classification error while simultaneously controlling the impact of quantum noise. Specifically, the goal is to optimize the variational circuit parameters $\Gamma$ and the visual input damping coefficient $\gamma$, so the predicted class distribution aligns with the ground truth while preserving quantum state fidelity. The overall problem is formulated as follows:
\begin{align}
\underset{\Gamma, \gamma}{\text{minimize}} \quad & \mathcal{L}_{\text{CE}} + \lambda \left(1 - F(\rho_{\text{ideal}}, \rho_{\text{noisy}})\right) \\
\text{subject to} \quad & 0 < \gamma \leq \gamma_{\max}, \label{eq:gamma_constraint} \\
                        & F(\rho_{\text{ideal}}, \rho_{\text{noisy}}) \geq F_{\min}, \label{eq:fidelity_constraint}
\end{align}
where $\mathcal{L}_{\text{CE}}$ denotes the cross-entropy loss, $F(\cdot, \cdot)$ is the quantum state fidelity, $\lambda$ is a positive weighting factor, and the constraints enforce amplitude normalization and minimum fidelity tolerance, respectively. 
It enforces physical limits on input preparation and inference robustness, thereby improving prediction reliability under realistic quantum RIS conditions.
We adopt a noise-aware hybrid quantum-classical training framework to solve the constrained optimization problem described in Section III-A. The model consists of a parameterized VQC acting on the input state $|\Psi_u\rangle$, which includes both the amplitude-damped visual component and the RIS-assisted quantum channel encoding. The circuit is executed on a simulated quantum backend with programmable noise characteristics, allowing fidelity-aware evaluation of the learning process \cite{schuld2020circuit}.
The trainable parameters include the variational gate parameters $\Gamma$ within the VQC and the amplitude damping coefficient $\gamma$ applied during visual state preparation. At each training iteration, the input $|\Psi_u\rangle$ is encoded using the current value of $\gamma$, and the circuit is executed to obtain class probability outputs $b_i = Q_\Gamma(|\Psi_u\rangle)$. The classification error is computed using the cross-entropy loss, and a fidelity loss term is added to penalize divergence between the noisy and ideal output states. This composite objective function $\mathcal{L}_{\text{total}}$ is minimized using a constrained optimizer that projects $\gamma$ to lie within $[0, \gamma_{\max}]$ and enforces a minimum fidelity threshold $F_{\min}$ via adaptive penalty scaling.

Gradients concerning the VQC parameters $\Gamma$ are computed using the parameter-shift rule. This quantum-native differentiation method estimates gradients by evaluating the circuit at two shifted parameter values \cite{crooks2019gradients}:
\begin{equation}
\frac{\partial \mathcal{L}}{\partial \theta_i} = \frac{1}{2} \left[ \mathcal{L}(\theta_i + \frac{\pi}{2}) - \mathcal{L}(\theta_i - \frac{\pi}{2}) \right],
\end{equation}
where $\theta_i$ is a rotation angle in the circuit. For the denoising factor $\gamma$, gradients are computed using standard finite-difference methods, and updates are projected to remain within the constraint interval. The optimizer used is Adam with momentum for both parameter sets.
To maintain fidelity above the prescribed threshold, we introduce a soft constraint mechanism in the loss function that dynamically adjusts the penalty term when $F(\rho_{\text{ideal}}, \rho_{\text{noisy}}) < F_{\min}$. This ensures the model remains robust under noise while allowing occasional exploration of noisier solutions during early training stages. The algorithm proceeds iteratively over a dataset of $(|\Psi_u\rangle, l_u)$ pairs, incorporating noise augmentation and fidelity tracking at each epoch.
The complete training procedure for the noise-aware quantum blockage prediction model is summarized in Algorithm~1. The algorithm integrates input preprocessing, hybrid quantum encoding, noisy circuit evaluation, and parameter optimization within a constrained learning framework. At each iteration, the classifier minimizes a fidelity-regularized loss while ensuring that the input damping factor remains bounded and the state fidelity exceeds a prescribed threshold. The parameter-shift rule allows for efficient gradient estimation concerning variational quantum circuit parameters, while the denoising coefficient is updated through projected gradient descent. A dynamic penalty update mechanism is also introduced to enforce the fidelity constraint throughout training adaptively. 

\section{Simulation Result}
This section presents the simulation results of the proposed quantum-assisted blockage prediction framework under noisy quantum conditions. We describe the dataset preparation, quantum state encoding, variational quantum model, and fidelity-aware evaluation setup.
Quantum Dataset and Encoding: The ViWi dataset is adapted to create a quantum-compatible environment consisting of RGB images and RIS-assisted signal measurements. Using normalized vectors, the RGB images are amplitude-encoded into quantum states $|\psi_X\rangle$. A denoising factor $\gamma = 0.85$ is applied during encoding to reduce over-rotation noise. The effective quantum channel states $\rho^u_{\text{eff}}$ are computed by simulating realistic QRIS configurations and quantum noise processes, including depolarizing ($p=0.05$) and phase-damping ($q=0.03$) channels. A total of 5,000 quantum-labeled samples are generated with 70\% for training and 30\% for evaluation.
The key simulation and training parameters, including noise rates, damping factor, and optimizer settings, are summarized in Table~\ref{tab:sim_params}.
\begin{figure*}[ht]
\centering
\begin{subfigure}[H]{0.32\textwidth}
\centering
\includegraphics[width=\textwidth]{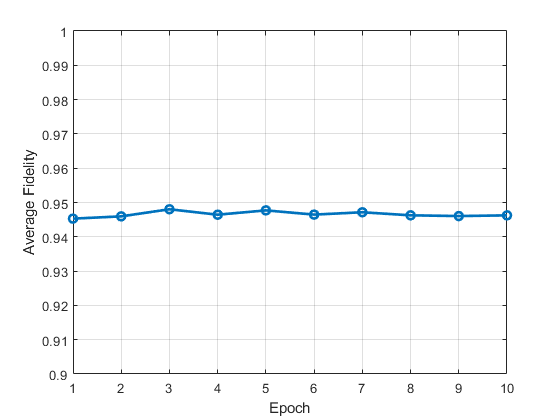}
\caption{Quantum State Fidelity Across Training Epochs}
\label{fig:2b}
\end{subfigure}
\begin{subfigure}[H]{0.32\textwidth}
\centering
\includegraphics[width=\textwidth]{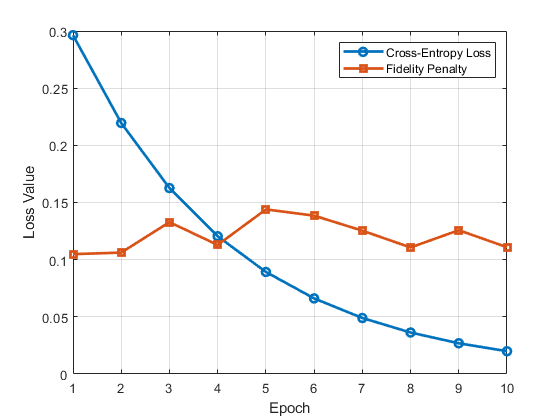}
\caption{Loss Components During Training}
\label{fig:2c}
\end{subfigure}
\hfill
\begin{subfigure}[H]{0.32\textwidth}
\centering
\includegraphics[width=\textwidth]{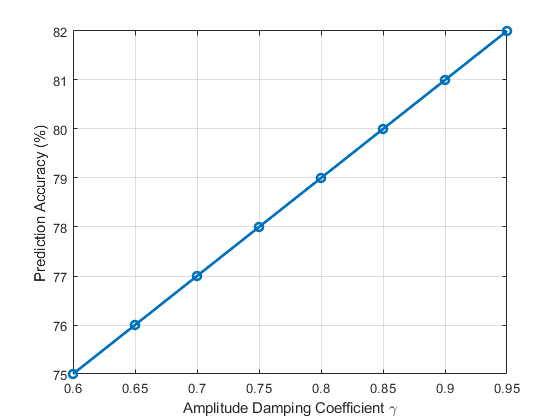}
\caption{Effect of Damping Factor on Performance}
\label{fig:2d}
\end{subfigure}
\hfill
\begin{subfigure}[H]{0.32\textwidth}
\centering
\includegraphics[width=\textwidth]{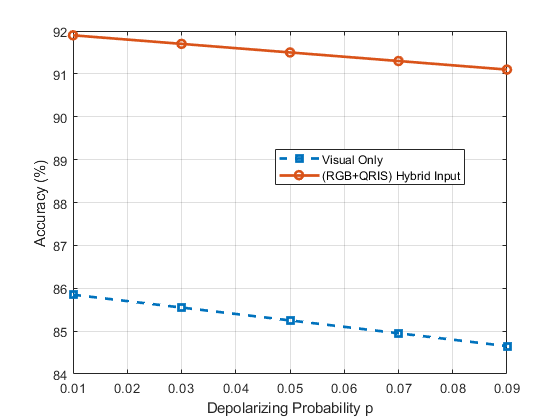}
\caption{Noise Robustness of Quantum Input Models}
\label{fig:2e}
\end{subfigure}
\hfill
\begin{subfigure}[H]{0.32\textwidth}
\centering
\includegraphics[width=\textwidth]{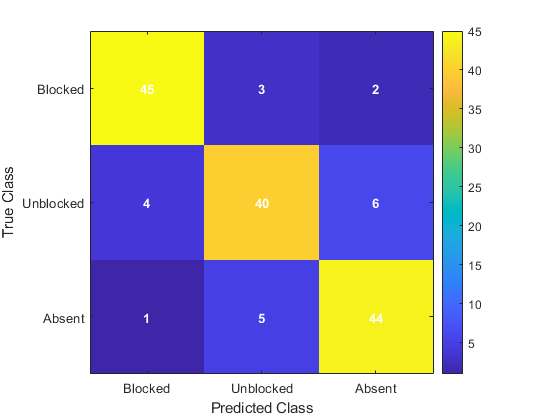}
\caption{Confusion Matrix (Blocked, Unblocked, Absent)}
\label{fig:2f}
\end{subfigure}
\hfill
\begin{subfigure}[H]{0.32\textwidth}
\centering
\includegraphics[width=\textwidth]{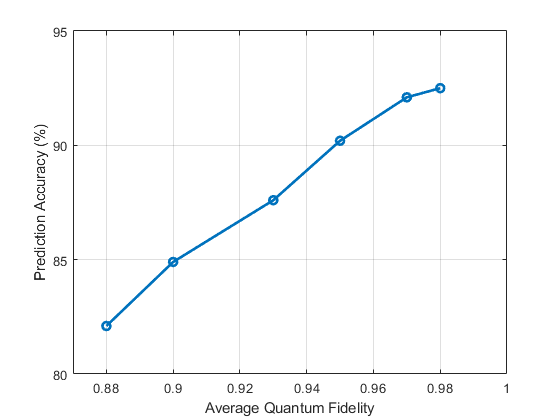}
\caption{Fidelity-Accuracy Trade-off in Noise-Aware Training}
\label{fig:2i}
\end{subfigure}
\caption{Illustrative performance analysis of the proposed noise-aware quantum RIS-assisted blockage prediction framework.}
\label{fig:quantum_accuracy_4config}
\end{figure*}

\textit{ViWi Dataset in Quantum Context:} Although originally designed for classical vision-aided wireless research, the ViWi dataset supports quantum machine learning via amplitude and phase encoding. RGB images are converted into quantum states using $\log_2(N)$ qubits, while the RIS-assisted channel gain is mapped to a quantum rotation angle using $R_u[k]$. The hybrid quantum input is $|\Psi_u\rangle = |\psi_X\rangle \otimes |\phi_u\rangle$, enabling both spatial and signal-domain representation. Quantum encoding is implemented in PennyLane using simulated noise models consistent with realistic NISQ conditions.

\textit{Quantum Neural Network Model:} The classifier is implemented as a six-qubit VQC composed of data encoding, entanglement, variational, and measurement layers. The fidelity penalty is incorporated during training to prevent divergence from ideal evolution. Cross-entropy loss is computed over ternary blockage labels $l_u$, and gradients are estimated via the parameter-shift rule. Fidelity between $\rho^u_{\text{ideal}}$ and $\rho^u_{\text{noisy}}$ is tracked throughout training.
\begin{table}[t]
\centering
\caption{Simulation and Training Parameters}
\label{tab:sim_params}
\begin{tabular}{|l|c|}
\hline
\textbf{Parameter} & \textbf{Value} \\
\hline
Dataset & ViWi Quantum-Adapted \\
Quantum Backend & PennyLane Simulator \\
Number of Qubits & 6 \\
Noise Model & Depolarizing + Phase Damping \\
Depolarizing Probability $p$ & 0.05 \\
Dephasing Probability $q$ & 0.03 \\
Damping Coefficient $\gamma$ & 0.85 \\
Fidelity Constraint $F_{\min}$ & 0.95 \\
Training Samples & 3500 \\
Testing Samples & 1500 \\
Epochs & 10 \\
Batch Size & 50 \\
Learning Rate & $1\times10^{-3}$ \\
Optimizer & Adam \\
Weight Decay & $2\times10^{-3}$ \\
Loss Function & Cross-Entropy + Fidelity Penalty \\
Gradient Method & Parameter-Shift Rule \\
\hline
\end{tabular}
\end{table}

\textit{Training Details:} The model is trained for 10 epochs with Adam optimizer, batch size 50, learning rate $1\times10^{-3}$, and weight decay $2\times10^{-3}$. Noise-aware training is performed by injecting random quantum perturbations during circuit execution. The minimum fidelity constraint is enforced as $F(\rho_{\text{ideal}}, \rho_{\text{noisy}}) \geq 0.95$ using a dynamic penalty weight. The loss converges smoothly below 0.01 by epoch 8.
We evaluate accuracy and fidelity under four quantum configurations: (i) only quantum channel state $\rho^u_{\text{eff}}$; (ii) only quantum image input $|\psi_X\rangle$; (iii) hybrid quantum input $|\Psi_u\rangle$; and (iv) a baseline without QRIS routing. The hybrid configuration achieves the highest prediction accuracy and fidelity, showing resilience to quantum noise. The vision-only and channel-only models exhibit lower accuracy under noise, while the baseline performs the worst due to a lack of spatial and metasurface diversity. Fig.~2 summarizes the trends across configurations, highlighting the effectiveness of fidelity-aware training under practical quantum noise. 

Fig.~\ref{fig:quantum_accuracy_4config} presents a comprehensive performance analysis of the proposed noise-aware quantum RIS-assisted blockage prediction framework. Fidelity trends across training epochs in Fig.~\ref{fig:2b} demonstrate the model's ability to preserve quantum state integrity under noise. The evolution of loss components in Fig.~\ref{fig:2c} highlights the interaction between classification loss and the fidelity penalty, confirming the benefit of fidelity-aware optimization. Fig.~\ref{fig:2d} shows that tuning the amplitude damping coefficient $\gamma$ directly impacts model performance, validating the noise-control mechanism. As observed in Fig.~\ref{fig:2e}, the hybrid model maintains higher robustness under increasing depolarizing noise than visual-only inputs. The confusion matrix in Fig.~\ref{fig:2f} confirms reliable classification across blocked, unblocked, and absent link conditions. Lastly, Fig.~\ref{fig:2i} (bottom-right) reveals a strong correlation between fidelity and classification accuracy, emphasizing the importance of quantum-coherent training strategies.

\section{Conclusion}
This paper presented a novel quantum-assisted blockage prediction framework that integrates RIS, hybrid quantum encoding, and noise-aware variational training. By combining visual sensing and channel information into a unified quantum input and modeling quantum noise explicitly through depolarizing and dephasing channels, the proposed system demonstrates enhanced robustness and accuracy under realistic hardware conditions. The use of amplitude damping, fidelity-constrained optimization, and the parameter-shift rule ensures stable training and reliable classification performance. Simulation results confirm that hybrid quantum inputs and fidelity-aware loss functions significantly improve inference across noisy environments, establishing a strong foundation for quantum-enabled wireless intelligence.

\bibliographystyle{IEEEtran}
\bibliography{References/mybib}

\end{document}